\documentclass[12pt]{amsproc} %%
\usepackage{amssymb,amsmath}
\usepackage{mathrsfs}

\newtheorem{definition}{Definition}[section]
\newtheorem{theorem}{Theorem}[section]

\newtheorem{cor}{Corollary}[section]
\newtheorem{remark}{Remark}[section]
\newtheorem{example}{Example}[section]
\newtheorem{step}{Step}
\def\Proof{{\bf Proof.} \begin{em} }
\def\endProof{\end{em}\hfill $\blacksquare$\par}

\def\bbk{\mathbb{K}}
\def\bbl{\mathbb{L}}

\def\scr#1{\mathscr{#1}}

\def\ca{\mathcal{A}}
\def\cb{\mathcal{B}}
\def\ci{\mathcal{I}}
\def\cf{\mathcal{F}}

\def\borel{{\mathcal{B}_+}}

\def\bigcups{\bigcup}

\def\<{\langle}
\def\>{\rangle}

\def\Cichon{Cicho{\'n}}
\def\Zeberski{{\.Z}eberski}
\def\Ralowski{Ra{\l}owski}

\begin{document}
\title{Complete nonmeasurability in regular families}
\author{Robert Ra{\l}owski\and Szymon \Zeberski}

\address{Robert \Ralowski \and Szymon \Zeberski, Institute of Mathematics, Wroc{\l}aw University of Technology, Wybrze\.ze Wyspia\'n\-skie\-go 27, 50-370 Wroc{\l}aw, Poland.}
\email{robert.ralowski@pwr.wroc.pl}

\email{szymon.zeberski@pwr.wroc.pl}

\thanks{ Authors would like to thank  prof. Jacek {\Cichon} 
        for many helpful suggestions.}
%\date{\wersja} %%{August 2004}

\begin{abstract}
We show that for a  $\sigma $-ideal $\ci$ with a
Borel base of subsets of an uncountable Polish space, if
$\ca$ is (in several senses) a "regular" family of subsets from $\ci $ then
there is a subfamily of $\ca$ whose union is completely
nonmeasurable i.e. its intersection with every Borel set not in $\ci $ does not
belong to the smallest $\sigma $-algebra containing all Borel sets and $\ci.$ Our results 
generalize   results from \cite{fourpoles} and \cite{fivepoles}.
\end{abstract}

\maketitle

\section{Notation and Terminology}

Throughout this paper, $X$, $Y$ will denote uncountable Polish spaces and $\cb(X)$ the Borel $\sigma $-algebra of $X.$  
We say that the ideal $\ci$ on $X$ has {\it Borel base}  if every  element $A\in\ci$ is contained in a Borel set in $\ci.$ (It is assumed that an ideal is always proper.)
The ideal consisting of all countable subsets of $X$ will be denoted by $[X]^{\le\omega }$ and the ideal of all meager subsets of $X$ will be denoted by $\bbk.$
Let $\mu $ be a continous probability measure on $X.$ The ideal consisting of all $\mu $-null sets will be denoted by $\bbl_\mu.$ By the following well known result, $\bbl_\mu $ can be identified with the $\sigma $-ideal of Lebesgue null sets.

\begin{theorem}[\cite{SMS}, Theorem 3.4.23]
If $\mu $ is a continous probability on $\cb (X),$ then there is a Borel isomorphism $h:X\rightarrow [0,1]$ such that for every Borel subset $B$ of $[0,1],$ $\lambda(B)=\mu(h^{-1}(B)),$ where $\lambda $ is a Lebesgue measure. 
\end{theorem}
 
\begin{definition}
We say that $(Z,\ci)$ is Polish ideal space if $Z$ is Polish uncountable space and $\ci$ is a $\sigma $-ideal on $Z$ having Borel base and containing all singletons. In this case, we set
$$\borel (Z)=\cb (Z)\setminus\ci.$$
A subset of $Z$ not in $\ci $ will be called a $\ci $-positive set; sets in $\ci $ will also be called $\ci $-null. Also, the $\sigma $-algebra generated by $\cb (Z)\cup\ci$ will be denoted by $\overline{\cb}(Z),$ called the $\ci $-completion of $\cb(Z).$
\end{definition}

It is easy to check that $A\in\overline{\cb }(Z)$ if and only if there is an $I\in\ci $ such that $A\bigtriangleup I$ (the symetric difference) is Borel.

\begin{example}
Let $\mu $ be a continous probability measure on $X.$ Then $(X,[X]^{\le\omega})$,
$(X,\bbk )$, $(X,\bbl_\mu )$ are Polish ideal spaces. 
\end{example}

\begin{definition}
A Polish ideal group is 3-tuple $(G,\ci,+)$ where $(G,\ci )$ is Polish ideal space and $(G,+)$ is an abelian topological group with respect to the Polish topology of $G.$
\end{definition}

\begin{definition}
Let $(X,\ci )$ be a Polish ideal space and $A\subseteq X$. We say that $A$ is $\ci$--nonmeasurable, if $A\notin\overline{\cb }(X).$ Further, we say that $A$ is completely $\ci $--nonmeasurable if
$$\forall B\in\borel (X) \;\;A\cap B\ne\emptyset\land A^c\cap B\ne\emptyset .
$$
\end{definition}

Clearly every completely $\ci $--nonmeasurable set is $\ci $--nonmeasurable. In the literature, completely $[X]^{\le\omega }$--nonmeasurable sets are called  Bernstein sets.
Also, note that  $A$ is completely $\bbl_\mu$--nonmeasurable if and only if the inner measure of $A$ is zero and the outer measure one.

For any set $E$, $|E|$ will denote the cardinality of $E.$

Let $(X,\ci)$ be a Polish ideal space and $\cf\subseteq\ci.$ We set 
$$\begin{array}{l@{\ =\,{} }l}
add(\ci) &  \min\{ |\ca|:\;\ca\subseteq \ci\land \bigcup \ca\notin \ci \}\\
cov(\ci) &  \min\{ |\ca|:\;\ca\subseteq \ci\land \bigcup \ca=X \}\\
cov(\cf) &  \min\{ |\ca|:\;\; \ca\subseteq\cf\land \bigcup \ca=X\}\\
cov_h(\ci) &  \min\{ |\ca|:\;\ca\subseteq \ci\land \exists B\in \borel(X) B\subseteq  
       \bigcup \ca \}\\
cov_h(\cf) &  \min\{ |\ca|:\;\ca\subseteq \cf\land \exists B\in \borel(X) B\subseteq 
       \bigcup \ca \}\\
\end{array}$$

An ideal $\ci $ is c.c.c. if every family of pairwise disjoint non-empty $\ci $-positive Borel sets  is countable.
Now let $(X,\ci )$ be a Polish ideal space with $\ci $ c.c.c. and $A\subseteq X.$
Let $\mathcal{A}$ be a maximal family of pairwise disjoint $\ci $-positive Borel sets
contained in $A^c.$ Set $B=(\bigcup\mathcal{A})^c.$ Then $B$ is Borel, $A\subseteq B$ and for every Borel set $C\supseteq A$, $B\setminus C\in\ci.$ Any such set $B$ is called a {\it Borel envelope} of $A$ and will be denoted by $[A]_\ci.$ Note that a Borel envelope of $A$ is unique modulo $\ci $ and it is minimal (modulo $\ci $) Borel set containing $A.$

It follows that $\overline{\cb }(X)$ is Marczewski complete (see \cite{SMS}, p.114). Therefore, it is closed under Souslin operation (see \cite{SMS}, Theorem 3.5.22). It follows that if $\ci $ is also c.c.c., $\overline{\cb }(X)$ contains all analytic sets.

For any set $F\subseteq X\times Y$ and $x\in X$, $y\in Y$ let 
$$F_x=\{ y\in Y:\; (x,y)\in F\}$$ 
and 
$$F^y=\{ x\in X:\; (x,y)\in F\}.$$
Further, for any $T\subseteq Y,$ we set
$$ F^{-1}(T)=\{x\in X: F_x\cap T\neq\emptyset\}.$$
A multifunction $F:X\rightarrow Y$ is called $\ca $--{\it measurable} if for every open set $U$ in $Y$, $F^{-1}(U)\in\ca,$ where $\ca $ is a $\sigma $-algebra on $X.$ 

Let $\pi $ be a partition of $X$ and $A\subseteq X.$  The smallest $\pi $-invariant subset of $X$ containing $A$ is called the {\it saturation} of $A$ and is denoted by $A^*.$ Thus,
$$A^*=\bigcup\{E\in\pi : E\cap A\neq\emptyset\}.$$
We call $\pi$ {\it Borel measurable} if the saturation of every open set is Borel; it is {\it strongly Borel measurable} if the saturation of every closed set is Borel measurable. Since $X$ is second countable, every strongly Borel measurable partition is Borel measurable.

The rest of our notations and terminology are standard. For other notation and terminology in Descriptive Set Theory we follow \cite{SMS}.

%%%%%%%%%%%%%%%%%%%%%%%%%%%%%%%%%%%%%%%%%%%%%%%%%
\section{Main results }
%%%%%%%%%%%%%%%%%%%%%%%%%%%%%%%%%%%%%%%%%%%%%%%%
The following results are the main results of the paper.

\begin{theorem}\label{selectors} Let $(X,\ci )$ be a Polish ideal space such that every set in $\cb_+(X)$ contains a $\ci $-positive closed set. Suppose  $\ca $ is a strongly Borel measurable partition of $X$ into $\ci $-null closed sets.
Then there is a subfamily $\ca_0\subseteq\ca$ such that $\bigcup\ca_0$ is completely $\ci $--nonmeasurable.
\end{theorem}

\begin{theorem}\label{open}
 Let $(X,\ci )$ be a Polish ideal space. Suppose $f:X\rightarrow Y$ is a $\overline{\cb }(X)$-measurable map such that for every $y\in Y$, $f^{-1}(y)\in\ci.$
Then there is a  $T\subseteq Y$ such that $f^{-1}(T)$ is completely
$\ci$--nonmeasurable.
\end{theorem}

\begin{theorem}\label{openccc}
 Let $(X,\ci )$ be a Polish ideal space with $\ci $ c.c.c. Let $F: X\rightarrow Y$ be a $\overline{\cb }(X)$-measurable multifunction such that for every $x\in X$, $F(x)$ is finite.
Then there exists a $T\subseteq Y$ such that $F^{-1}(T)$ is completely
$\ci$--nonmeasurable.
\end{theorem}

\begin{theorem}\label{finccc}
 Let $(X,\ci )$ be a Polish ideal space with $\ci $ c.c.c. Suppose $F$ is an analytic subset
of $X\times Y$  satisfying the following conditions:
\begin{enumerate}
\item $(\forall y\in Y)(F^y\in\ci )$;
\item $X\setminus \pi_X(F)\in\ci$, where $\pi_X:X\times Y\rightarrow X$ is the projection 
      map;
\item $(\forall x\in X)( |F_x|<\omega ).$
\end{enumerate}
Then there exists a $T\subseteq Y$ such that $F^{-1}(T)$ is completely $\ci$--nonmeasurable.
\end{theorem}

These results generalize results from \cite{fourpoles} and \cite{fivepoles}. In the next section, we present the proofs of our theorems.

%%%%%%%%%%%%%%%%%%%%%%%%%%%%%%%%%%%%%%%%%%%%
\section{Proofs of the main results}
%%%%%%%%%%%%%%%%%%%%%%%%%%%%%%%%%%%%%%%%%%%%

One of the key ideas of this paper is the following theorem (see \cite{fivepoles}).
For reader's convenience we will give the proof of it.

\begin{theorem}\label{klucz} Let $(X,\ci )$ be a Polish ideal space. Assume that a family $\ca\subseteq\ci $ satisfies the following conditions:
\begin{enumerate}
\item $X\setminus\bigcup\ca\in\ci$,
\item $Z=\{ x\in X:\;\;\bigcup\{ A\in\ca:\; x\in A\}\notin\ci\}\in\ci$,
\item $cov_h(\cf)=2^\omega$, where $\cf=\{ \bigcups\{ A\in\ca:\;\; x\in A\} :\;\; x\in 
       X\setminus Z\}$.
\end{enumerate}
Then there exists a subfamily $\ca_0\subseteq\ca $ such that $\bigcup\ca_0$ is completely $\ci $--nonmeasurable.
\end{theorem}
\proof 
First of all, we can assume that $Z=\emptyset$ in the second assumption. Now, let us enumerate the family of all positive Borel sets with respect to the ideal $\ci $ i.e. $\borel (X)=\{ B_\alpha:\;\;\alpha<2^\omega\}.$
By transfinite induction we will construct a sequence
$$
\< (d_\xi,A_\xi)\in B_\xi\times \ca:\;\;\xi<2^\omega \>
$$
satisfying the following conditions
\begin{enumerate}
\item $ A_\xi\cap B_\xi\neq\emptyset $,
\item $ d_\xi\notin\bigcups_{\alpha <2^\omega } A_\alpha $.
\end{enumerate}
Assume that we have constructed a sequence $\< (d_\xi,A_\xi)\in B_\xi\times \ca:\;\;\xi<\alpha \>.$
Since $\bigcups_{\xi<\alpha}\{ A\in \scr{A}:\; d_\xi\in A\}$ does not cover any positive Borel set, we are able to find $a_\alpha\in B_\alpha\setminus \bigcups_{\xi<\alpha}\{ A\in \ca:\; d_\xi\in A\}.$ Let $A_\alpha $ be any element of $\ca $ such that $a_\alpha\in A_\alpha$ and find $d_\alpha\in B_\alpha\setminus \bigcup_{\xi\le\alpha} A_\xi$. 
It finishes $\alpha$ step of our construction. 

Now, let us define $\ca_0=\{A_\xi:\;\;\xi\in 2^\omega\}.$ For every positive Borel set we have that $\bigcup\ca_0\cap B\neq\emptyset $ and $\{d_\xi : \xi\in 2^\omega\}\cap B\neq\emptyset $. Moreover,  $\{d_\xi:\;\;\xi\in 2^\omega\}\cap\bigcup\ca_0=\emptyset.$
It shows that $\bigcup\ca_0$ is completely $\ci $--nonmeasurable. 
\endproof

\begin{remark} We can replace the last assumption in Theorem \ref{klucz}  by the set theoretic assumption $cov_h(\ci)=2^\omega$.
\end{remark}

As a corollary we have:
\begin{cor}[ZFC+CH] Let $(X, \ci )$ be a Polish ideal space. Let $\ca\subseteq \ci$ be a  point-countable  family i.e. $\forall x\in X\; |\{ A\in\ca:\;\;x\in A\}|\le\omega$ and $\bigcup \ca =X$.  Then there exists a subfamily $\ca_0\subseteq\ca $ such that $\bigcup\ca_0$ is completely $\ci $--nonmeasurable.
\end{cor}

It is also known that above corollary is independent from $ZFC$ theory (see \cite{fremlin}).

\proof[Proof of Theorem \ref{selectors}]
By Theorem \ref{klucz}, it is sufficient to prove that 
\linebreak
$cov_h(\ca )=2^\omega.$ Towards proving this, take any $B\in\cb_+(X).$ Let$F\subseteq B$ be a $\ci $-positive closed set. Let
$$\pi=\{E\cap F: E\in\cf\}.$$
Note that $\pi $ is uncountable and strongly Borel measurable partition of $F$ into closed sets. Since every strongly Borel measurable partition is Borel measurable, it is Borel measurable. Hence, it admits a Borel cross-selection $S$ (see \cite{SMS}, Theorem 5.4.3, see \cite{onepole}). Clearly $S$ is uncountable and, therefore of cardinality $2^\omega.$ This implies that $|\pi |=2^\omega.$ 
\endproof

As a corollary we get the following result for Polish groups:
\begin{cor} Let $(G, \ci, +)$ be a compact Polish ideal group. Suppose $\ci $ is closed under translations. Assume that
each set from $\cb_+(G)$ contains a $\ci $-positive closed set. Let $H<G$ be a perfect subgroup and $H\in\ci$. 
Then there exists a $T\subseteq G$ such that $T+H$ is completely $\ci$--nonmeasurable in $G$.
\end{cor}
\proof 
This follows from Theorem \ref{selectors} by taking $\ca $ to be the set of all left cosets of $H.$
\endproof

To prove Theorem \ref{open}, we need the following result from \cite{fourpoles}.

\begin{theorem}[Brzuchowski, \Cichon, Grzegorek, Ryll-Nardzewski]
 Let $(X,\ci)$ be a Polish ideal space and $\ca\subseteq\ci $ a point-finite cover of $X.$ 
Then there is a subfamily $\ca_0\subseteq\ca $ whose union is not in $\overline{\cb }(X).$
\end{theorem}

\proof[Proof of Theorem \ref{open}]
Fix a countable base $\{U_n\}$ for the topology of $Y.$ For each $n,$ let $I_n\in\ci $ such that $f^{-1}(U_n)\bigtriangleup I_n$ is Borel. Let $X'=X\setminus\bigcup_{n}I_n.$ Then $f:X'\rightarrow Y$ is Borel. Thus, without any loss of generality, we assume that $f$ is Borel measurable.

Now, let $B\in\cb_+(X).$ Set
$$A=\pi_Y((B\times Y)\cap graph(f)).$$ 
Then $A,$ being analytic, is either countable or of cardinality $2^\omega.$
If $A$ were countable, $B$ is covered by countable subfamily of $\ci,$ a contradiction, Thus,
$cov_h\{ f^{-1}(y):y\in Y\}=2^\omega.$ Our result now follows from Theorem \ref{klucz}.
\endproof

\begin{theorem}\label{analityczne} Let $(X,\ci )$ be a Polish ideal space. Let $I\in\ci$ and $f:X\setminus I\rightarrow Y$ a Borel map such that for every $y\in Y$, $f^{-1}(y)$ is $\ci $-null. Then  there is a $T\subseteq Y$ such that $f^{-1}(T)$ is completely $\ci$--nonmeasurable set.
\end{theorem}
\proof
Let $B\supseteq I$ be a Borel $\ci$-null set. Now apply Theorem \ref{open} to $f\lceil(X\setminus B).$
\endproof

The next theorem is a technical result which helps us to prove stronger theorems in case $\ci $ is c.c.c.
\begin{theorem}\label{analityczneccc}
Let $(X,\ci )$ be a Polish ideal space with $\ci $ c.c.c. Assume that we have a family $\cf\subseteq\ci $ satisfying the following conditions:
\begin{enumerate}
\item $\cf$ is point-finite;
\item $(\forall B\in\cb_+(X))(B\subseteq [\bigcup\cf ]_\ci \rightarrow |\{F\in\cf : F\cap B\neq\emptyset \}|=2^\omega ).$  
\end{enumerate}
Then there exists a subfamily $\cf'\subseteq\cf $ such that $\bigcup\cf' $ is completely $\ci $--nonmeasurable in $[\bigcup\cf ]_\ci.$
\end{theorem}
\proof
\begin{step}
There exists a subfamily $\cf_0\subseteq\cf $ having the following properties
\begin{enumerate}
\item $[\bigcup\cf_0 ]_\ci =[\bigcup\cf ]_\ci, $
\item $(\forall B\in\cb_+(X))(B\subseteq\bigcup\cf_0 \rightarrow\min\{|\ca | :    
       \ca\subseteq\cf_0\wedge B\subseteq\bigcup\ca\}=2^\omega ).$   
\end{enumerate}
\end{step}
\proof
Let us recall that for a set $D\subseteq X$ a symbol $]D[_\ci$ denotes a maximal Borel set (mod $\ci $) contained in $D$. 
We will construct a sequence $(\ca_n)$ satistying the following conditions
\begin{enumerate}
\item $|\ca_n|<2^\omega,$
\item $\ca_n\subseteq\cf\setminus\bigcup_{i<n}\ca_i,$
\item $]\bigcup\ca_n[_\ci$ is maximal element in the family $\{]\bigcup\ca [_\ci : |\ca |<2^\omega
       \wedge \ca\subseteq\cf\setminus \bigcup_{i<n}\ca_i\}.$
\end{enumerate}
Notice that the existance of the maximal element in the family $\{]\bigcup\ca [_\ci : |\ca |<2^\omega \wedge \ca\subseteq\cf\setminus \bigcup_{i<n}\ca_i\}$ is implied by the c.c.c property of the ideal $\ci.$

We finish the construction if $\{]\bigcup\ca [_\ci : |\ca |<2^\omega
       \wedge \ca\subseteq\cf\setminus \bigcup_{i<n}\ca_i\}=\{\emptyset\}.$
Our construction has to end up after finitely many steps. Notice that $]\bigcup\ca_{n+1} [_\ci\subseteq ]\bigcup\ca_{n} [_\ci $ and $]\bigcup\ca_{n} [_\ci\neq\emptyset.$ So, assuming that there is infinitely many $\ca_n$'s we find a point $x\in X$ which belongs to infinitely many $\bigcup\ca_n$'s. Then $x$ belongs to infinitely many members of $\cf,$ what gives a contradiction with point-finiteness of the family $\cf.$ So, our construction ends up after $k$ steps ($k<\omega$).

Now, put $\cf_0=\cf\setminus\bigcup\{\ca_n : n\le k\}.$ It is a desired family.         
\endproof

\begin{step}
There exists a subfamily $\cf'\subseteq\cf_0 $ such that $\bigcup\cf'$ is completely $\ci$-non\-me\-a\-surable in $[\bigcup\cf_0]_\ci.$ 
\end{step}
\proof
Let us enumerate two families of positive Borel sets. Namely,
$$\cb^0=\{B^0_\alpha : \alpha < 2^\omega\}=\left\{B\in\cb_+(X) : B\subseteq \left[\bigcup\cf_0\right]_\ci\setminus \left]\bigcup\cf_0\right[_\ci\right\},$$
$$\cb^1=\{B^1_\alpha : \alpha < 2^\omega\}=\left\{B\in\cb_+(X) : B\subseteq  \left]\bigcup\cf_0\right[_\ci\right\}.$$
By transfinite induction we construct a sequence
$$
( (F^0_\xi,F^1_\xi, d_\xi)\in \cf_0\times\cf_0\times B^1_\xi:\;\;\xi<2^\omega )
$$
satisfying the following conditions
\begin{enumerate}
\item $F^0_\xi\cap B^0_\xi\neq\emptyset, \quad F^1_\xi\cap B^1_\xi\neq\emptyset,$
\item $d_\xi\notin\bigcup_{\xi<2^\omega } (F^0_\xi\cup F^1_\xi ).$
\end{enumerate}
Assume that we have constructed a sequence $( (F^0_\xi,F^1_\xi, d_\xi)\in \cf_0\times\cf_0\times B^1_\xi:\;\;\xi<\alpha )$
Since $|\{ F\in \cf_0:\; d_\xi\in F \mbox{ for some }  \xi<\alpha\}|<2^\omega$, we are able to find $F^0_\alpha,  F^1_\alpha $ such that $F^0_\alpha, F^1_\alpha\notin \{ F\in \cf_0:\; d_\xi\in F \mbox{ for some }  \xi<\alpha\}$ and $F^0_\alpha\cap B^0_\alpha\neq\emptyset, F^1_\alpha\cap B^1_\alpha\neq\emptyset.$  What is more $\bigcup\{F^0_\xi, F^1_\xi:\xi\le\alpha\}$ does not cover $B^1_\alpha.$ So, we can pick 
$d_\alpha\in B^1_\alpha\setminus \bigcup\{F^0_\xi, F^1_\xi:\xi\le\alpha\}.$ 
It finishes $\alpha$ step of our construction. 

Now, let us define $\cf'=\{F^0_\xi, F^1_\xi:\;\;\xi\in 2^\omega\}.$ We have that $\bigcup\cf'$ has not empty intersection with any positive Borel set contained in $[\bigcup\cf_0 ]_\ci $ and $\{d_\xi : \xi\in 2^\omega\}$ has not empty intersection with every positive Borel set contained in $]\bigcup\cf_0 [_\ci $. Moreover,
$\{d_\xi:\;\;\xi\in 2^\omega\}\cap\bigcup\cf'=\emptyset$ It implies that $\bigcup\cf'$ does not contain any positive Borel set.
It shows that $\bigcup\cf'$ is completely $\ci $--nonmeasurable in $[\bigcup\cf_0]_\ci.$  
\endproof
Since $[\bigcup\cf ]_\ci=[\bigcup\cf_0]_\ci,$ it finishes the proof. 
\endproof

\begin{remark}
Assuming that $cov(\ci )>\omega_1 $ we can prove the same theorem for wider class of families. Namely, it is enough to assume that a family $\cf\subseteq\ci $ is point-countable, i.e. $(\forall x\in X)(|\{F\in\cf : x\in f\}|\le\omega.$ Since $cov(\ci )>\omega_1, $ there is a point which belongs to $\omega_1 $ many Borel sets with the same envelope.
\end{remark}

\proof[Proof of Theorem \ref{openccc}]
By an argument contained in the proof of Theorem \ref{open}, without loss of generality, we can assume that $F^{-1}(U)$ is Borel for every open set $U$ in $Y.$ Fix any $B\in\cb_+(X).$
By Kuratowski--Ryll-Nardzewski selection theorem (see \cite{SMS}, Theorem 5.2.1, see \cite{twopoles}), $F\lceil B$ admits a Borel selection $s.$ The range of $s,$ being uncountable, is of cardinality $2^\omega.$ This implies that the condition (2) of Theorem \ref{analityczneccc} is satisfied by $\cf=\{F^{-1}(y) : y\in Y\}.$ Since each $F(x)$ is finite, $\cf $ is point-finite. The result now follows from Theorem \ref{analityczneccc}. 
\endproof

\proof[Proof of Theorem \ref{finccc}]
Without loss of generality, we can assume that 
\linebreak
$\pi_X(F)=X.$ Since $I$ is c.c.c., every analytic set in $X$ is in $\overline{\cb }(X)$ (see Section 1). It follows that $F$ is the graph of $\overline{\cb }(X)$-measurable, finite set valued multifunction. The result follows from Theorem \ref{openccc}.
\endproof

\vskip 1cm

\end{document}